\documentclass[11pt]{article}
\usepackage{euscript}

\usepackage{amsfonts}
\usepackage{amsmath}
\usepackage{amssymb}
\usepackage{latexsym,color}
\usepackage{amsbsy}
\usepackage{graphicx}

\newtheorem{thm}{Theorem}[section]
\newtheorem{con}{Conjecture}[section]
\newtheorem{lemma}[thm]{Lemma}
\newtheorem{cor}[thm]{Corollary}
\newtheorem{pro}[thm]{Proposition}
\newtheorem{example}[thm]{Example}
\newtheorem{definition}[thm]{Definition}

\newtheorem{remark}[thm]{Remark}
\newtheorem{Algorithm}[thm]{Algorithm}

\newcommand{\comment}[1]{}
\newcommand{\cIq}{{\mathcal I}_q}
\newcommand{\cNq}{{\mathcal N}_q}

\newcommand{\cS}{\mathcal S}
\newcommand{\cH}{\mathcal H}

\newcommand{\Aq}{A_q}
\newcommand{\Bq}{B_q}

\newcommand{\Cq}{C_q}

\newcommand{\Fq}{{\mathbb F}_q}
\newcommand{\C}{{\mathbb C}}

\newcommand{\pr}{\parallel}

\newcommand{\kq}[1]{{{ {[#1]_q} }}}
\newcommand{\qb}[2]{{{ {{#1}\brack {#2}}_q }}}

\newcommand{\oc}[1]{{{ \overline{#1} }}}

\newcommand{\ncom}{\newcommand}
\ncom{\ns}{\normalsize}
\ncom{\la}{\lambda}
\ncom{\bm}{\boldmath}
\ncom{\noi}{\noindent}
\ncom{\bq}{\begin{equation}}
\ncom{\eq}{\end{equation}}
\ncom{\beqn}{\begin{eqnarray*}}
\ncom{\eeqn}{\end{eqnarray*}}
\ncom{\ba}{\begin{array}}
\ncom{\ul}{\underline}

\ncom{\ea}{\end{array}}
\ncom{\beq}{\begin{eqnarray}}
\ncom{\eeq}{\end{eqnarray}}
\ncom{\nno}{\nonumber}
\ncom{\hs}{\mbox{\hspace{.25cm}}}
\ncom{\rar}{\rightarrow}
\ncom{\Rar}{\Rightarrow}
\ncom{\noin}{\noindent}
\ncom{\bc}{\begin{center}}
\ncom{\ec}{\end{center}}
\ncom{\sz}{\scriptsize}
\ncom{\fpd}{\Phi(\pi^{'})}
\ncom{\fp}{\Phi(\pi) }
\ncom{\nk}{\left< \begin{array}{c}
                       n\\k \end{array} \right>}
\ncom{\nd}{1^{'},2^{'},\cdots,n^{'}}
\ncom{\R}{I\!\!R}
\ncom{\de}{\bigtriangleup (F_{2n},\leq)}
\ncom{\del}{\bigtriangleup}
\ncom{\cov}{<\!\!\!\!\cdot }
\ncom{\bt}{\begin{thm}}
\ncom{\bcon}{\begin{con}}
\ncom{\et}{\end{thm}}
\ncom{\econ}{\end{con}}
\ncom{\bl}{\begin{lemma}}
\ncom{\el}{\end{lemma}}
\ncom{\bco}{\begin{cor}}
\ncom{\ds}{\displaystyle}
\ncom{\eco}{\end{cor}}
\ncom{\bp}{\begin{pro}}
\ncom{\ep}{\end{pro}}
\ncom{\bex}{\begin{example}}
\ncom{\eex}{\end{example}}
\ncom{\bd}{\begin{definition}}
\ncom{\ed}{\end{definition}}
\ncom{\brm}{\begin{remark}}
\ncom{\erm}{\end{remark}}
\ncom{\bal}{\begin{Algorithm}}
\ncom{\eal}{\end{Algorithm}}
\ncom{\ol}{\overline}
\ncom{\wh}{\widehat}

\ncom{\pf}{\noi {\bf Proof  }}
\ncom{\be}{\begin{enumerate}}
\ncom{\ee}{\end{enumerate}}
\ncom{\s}{\subset}
\ncom{\T}{{\cal T}}
\ncom{\B}{{\cal B}}
\ncom{\A}{{\cal A}}
\ncom{\Z}{{\cal Z}}

\title{\Large{{{\bf The Goldman-Rota identity and the Grassmann scheme}}}}

\author{{ {\bf Murali K. Srinivasan}} \\
{\em  {Department of Mathematics}}\\
{\em  {Indian Institute of Technology, Bombay}}\\
{\em  {Powai, Mumbai 400076, INDIA}}\\
{\bf  \texttt{mks@math.iitb.ac.in}}\\
{\bf  \texttt{murali.k.srinivasan@gmail.com}}
}

\begin{document}
\date{}
\maketitle

\begin{abstract}

We inductively construct an explicit (common) orthogonal eigenbasis for the
elements of the Bose-Mesner algebra of the Grassmann scheme. The main step
is a constructive, linear algebraic interpretation of the
Goldman-Rota recurrence for the number of subspaces of a finite vector
space. This interpretation shows that the up operator on subspaces has
an explicitly given recursive structure. Using this we 
inductively construct an explicit orthogonal symmetric 
Jordan basis with respect to the
up operator and write down the singular values, i.e., the ratio of the
lengths of the successive vectors in the Jordan chains. 
The collection of all vectors in this basis of a fixed rank
forms a (common) orthogonal eigenbasis for the elements 
of the Bose-Mesner algebra of the Grassmann scheme. We also pose a bijective
proof problem on the spanning trees of the Grassmann graphs.

\end{abstract}

{{\bf  \section {  Introduction }}}  

This paper presents {\em constructive} and {\em explicit} proofs of two 
basic linear algebraic results on the subspace lattice.

The first result concerns the recursive structure of the 
up operator on subspaces. It is an elementary observation that
the up operator (or equivalently, incidence matrices) on subsets of a
$n+1$-set can be built from two copies of the up operator on subsets of a
$n$-set. The main purpose of this paper is to extend this inductive approach
to the subspace lattice. A classical identity of Goldman and Rota suggests
that the up operator on subspaces of a $n+1$-dimensional vector space over
$\Fq$ can be built from two copies of the up operator in dimension $n$ and
$q^n-1$ copies of the up operator in dimension $n-1$. Let us make this
precise.

Let $B(n)$ denote the collection of all subsets of the set $\{1,2,\ldots
,n\}$. Partially order $B(n)$ by inclusion (our poset terminology follows
{\bf\cite{st2}}). Then $|B(n)| = 2^n$ and the identity $2^{n+1}=2\cdot 2^n$
has the following poset theoretic interpretation: we can write $B(n+1)$ as a
disjoint union $B(n+1) = B(n) \cup [\{n+1\}, \{1,\ldots ,n+1\}]$, where
the interval $[\{n+1\}, \{1,\ldots ,n+1\}]$ is order isomorphic to $B(n)$.
 
Let $\Fq^n$ denote the $n$-dimensional vector space of all column
vectors of length $n$ over $\Fq$ and let $\Bq(n)$ denote the collection of
all subspaces of $\Fq^n$. Partially order $\Bq(n)$ by inclusion. The number
of subspaces in $\Bq(n)$ having dimension $k$ is the $q$-binomial
coefficient $\qb{n}{k}$ and the total number of subspaces is the {\em Galois
number}
$$G_q(n)=\sum_{k=0}^n \qb{n}{k}.$$

The Goldman-Rota identity {\bf\cite{gr,kc,nsw}} is the recursion
\beq \label{gri}
G_q(n+1) &=& 2G_q(n) + (q^n - 1)G_q(n-1),\;n\geq
1,\;\;\;G_q(0)=1,\;G_q(1)=2. \eeq
We identify $\Fq^n$ with 
the subspace of all vectors in $\Fq^{n+1}$ with last component zero.
Put $t=q^n - 1$.
Motivated by the $B(n)$ case we can ask for the
following poset theoretic interpretation of (\ref{gri}): is it
possible to write $\Bq(n+1)$ as a disjoint union
\beq \label{du}
\Bq(n+1) &=& \Bq(n) \cup S_0 \cup S_1 \cup \cdots \cup S_t ,
\eeq
where $S_0,\ldots ,S_t$ are intervals in $\Bq(n+1)$, with $S_0$ order
isomorphic to $\Bq(n)$ and $S_1,\ldots ,S_t$ order isomorphic to
$\Bq(n-1)$. At least for $q=2$ and $n\geq 4$ the answer is no, as shown in
{\bf\cite{hh}}.

We show that we can get a poset theoretic
interpretation of (\ref{gri}) by considering a linear analog of (\ref{du}).
Moreover, the linear analog of the decomposition (\ref{du}) can be
explicitly given.

Let $P$ be a finite {\em graded poset} with
{\em rank function}
$r: P\rar {\mathbb N}=\{0,1,2,\ldots \}$. The {\em rank} of $P$ is
$r(P)=\mbox{max}\{r(x): x\in P\}$ and,
for $i=0,1,\ldots ,r(P)$, $P_i$ denotes the set of elements of $P$ of rank
$i$. For a subset $S\subseteq P$, we set $\mbox{rankset}(S) = \{r(x):x\in 
S\}$.

For a finite set $S$, let $V(S)$ denote the complex vector space with $S$ as
basis. Let $v=\sum_{x\in S}\alpha_x\,x,\,\alpha_x\in \C$ be an element of
$V(S)$. By the {\em support} of $v$ we mean the subset $\{x\in S :
\alpha_x\not= 0\}$. 

Let $P$ be a graded poset with $n=r(P)$. Then we have
$ V(P)=V(P_0)\oplus V(P_1) \oplus \cdots \oplus V(P_n)$ (vector space direct
sum).
An element $v\in V(P)$ is {\em homogeneous} if $v\in V(P_i)$ for some $i$,
and if $v\not= 0$, we extend the notion of rank to nonzero homogeneous
elements by writing
$r(v)=i$. Given an element $v\in V(P)$, write $v=v_0 + \cdots +v_n,\;v_i \in
V(P_i),\;0\leq i \leq n$. We refer to the $v_i$ as the {\em homogeneous
components} of $v$. A subspace $W\subseteq V(P)$ is {\em homogeneous} if it
contains the homogeneous components of each of its elements. For a
homogeneous subspace $W\subseteq V(P)$ we set $\mbox{rankset}(W)=\{r(v) : v
\mbox{ is a nonzero homogeneous element of } W\}$.

The {\em
up operator}  $U:V(P)\rar V(P)$ is defined, for $x\in P$, by
$U(x)= \sum_{y} y$,
where the sum is over all $y$ covering $x$. We denote the up operator on
$V(\Bq(n))$ by $U_n$. For a finite vector space $X$ over $\Fq$ we denote by
$\Bq(X)$ the set of all subspaces of $X$ and we denote by $U_X$ the up
operator on $V(\Bq(X))$.

Let $\langle , \rangle$ denote the standard inner product on $V(P)$,
i.e.,
$\langle x,y \rangle =
\delta (x,y)$ (Kronecker delta), for $x,y\in P$.
The {\em length} $\sqrt{\langle v, v \rangle }$ of $v\in V(P)$ is denoted
$\pr v \pr$.

Let $(V,f)$ be a pair consisting of a finite dimensional vector space
$V$ (over $\C$) and a linear operator $f$ on $V$. Let $(W,g)$ be
another such pair. By an isomorphism of pairs $(V,f)$ and $(W,g)$ we mean a
linear isomorphism $\theta : V 
\rar W$ such that $\theta(f(v)) = g(\theta(v)),\;v\in V$.

We give $V(\Bq(n))$ (and $V(\Bq(X))$ for a finite vector space $X$ over
$\Fq$) the standard inner product.
In Section 2 we prove the following 
result on the recursive structure of the pair $(V(\Bq(n)),U_n)$. 
Taking dimensions we get (\ref{gri}).
\bt \label{grv}
Set $t=q^n - 1$. 
There is an explicit orthogonal direct sum decomposition
\beq \label{odsd}
V(\Bq(n+1)) &=& V(\Bq(n)) \oplus W(0) \oplus W(1) \oplus \cdots \oplus W(t),
\eeq
where

\noi (i) $W(0),\ldots ,W(t)$ are $U_{n+1}$-closed (i.e., closed under the
action of $U_{n+1}$) homogeneous subspaces of
$V(\Bq(n+1))$ with $\mbox{rankset}\,(W(0))=\{1,\ldots ,n+1\}$ and 
$\mbox{rankset}\,(W(i))=\{1,\ldots ,n\},$ for $i=1,\ldots ,t$.

\noi (ii) $V(\Bq(n)) \oplus W(0)$ is $U_{n+1}$-closed and
there is an explicit linear map $\theta_n : V(\Bq(n)) \rar W(0)$
that is an isomorphism of pairs 
$(V(\Bq(n)), q U_n)$ and 
$(W(0),U_{n+1})$, sending homogeneous elements to homogeneous elements,
increasing rank by one and satisfying
\beq \label{ti} 
U_{n+1}(v) &=& U_n(v) + \theta_n(v),\,\,v\in V(\Bq(n)), \\
\langle \theta_n(w), \theta_n(v) \rangle &=& 
q^{n-k} \langle w,v\rangle,\,\,w,v\in
V(\Bq(n)_k),\,0\leq k \leq n.
\eeq

\noi (iii) For $i=1,\ldots ,t$ there is an explicit linear map 
$\gamma_{n-1}(i): V(\Bq(n-1)) \rar W(i)$ that is an isomorphism of pairs
$(V(\Bq(n-1)), U_{n-1})$ and
$(W(i),U_{n+1})$, sending homogeneous elements to homogeneous elements,
increasing rank by one and satisfying
\beq  \label{ti2}
\langle \gamma_{n-1}(i)(w), \gamma_{n-1}(i)(v) \rangle &=& 
q^{n+k} \langle w,v\rangle,\,\,w,v\in
V(\Bq(n-1)_k),\,0\leq k \leq n-1.
\eeq

\et

Our second main result is concerned with explicit construction of orthogonal
symmetric Jordan bases. Let $P$ be a finite graded poset with rank function $r$.
A {\em graded Jordan chain} in $V(P)$ is a sequence
\beq \label{gjc}
&s=(v_1,\ldots ,v_h)&
\eeq 
of nonzero homogeneous elements of $V(P)$
such that $U(v_{i-1})=v_i$, for
$i=2,\ldots h$, and $U(v_h)=0$ (note that the
elements of this sequence are linearly independent, being nonzero and of
different ranks). We say that $s$ {\em
starts} at rank $r(v_1)$ and {\em ends} at rank $r(v_h)$.
A {\em graded Jordan basis} of $V(P)$ is a basis of $V(P)$
consisting of a disjoint union of graded Jordan chains    
in $V(P)$. 
The graded Jordan
chain (\ref{gjc}) is said to be a {\em symmetric Jordan chain} (SJC) if
the sum of the starting and ending ranks of $s$ equals $r(P)$, i.e.,   
$r(v_1) + r(v_h) = r(P)$ 
if $h\geq
2$, or $2r(v_1)= r(P)$ if $h=1$.
A {\em symmetric Jordan basis} (SJB) of $V(P)$ is a basis of $V(P)$
consisting of a disjoint union of symmetric Jordan chains 
in $V(P)$.

Using Theorem \ref{grv} we prove the following result in Section 3.

\bt \label{mt3}

There is an algorithm to inductively construct an explicit orthogonal
SJB $J_q(n)$ of $V(\Bq(n))$. When expressed in the standard basis the
vectors in $J_q(n)$ have coefficients that are integral multiples of
$q${\em th} roots of unity. In particular, the coefficients are integral
when $q=2$. 

Let $0\leq k \leq n/2 $ and let
$(x_k,\ldots ,x_{n-k})$ be any SJC
in $J_q(n)$ starting at rank $k$ and ending at rank $n-k$. Then
we have, for $k\leq u < n-k$,
\beq \label{sv}
\frac{\pr x_{u+1} \pr}{\pr x_u \pr} & = &
\sqrt{q^k\kq{u+1-k}\kq{n-k-u}}\;\,.
\eeq
\et 

A standard argument (recalled in Section 3) shows that the set $\{v\in
J_q(n) : r(v)=m\}$ forms a common  orthogonal eigenbasis for the elements
of the Bose-Mesner algebra
of the Grassmann scheme of $m$-subspaces.

The numbers on the right hand side of (\ref{sv}) are called the {\em
singular values} of the up operator. These are important for applications.
The {\em existence} of an orthogonal SJB of $V(\Bq(n))$ 
satisfying (\ref{sv}) was
first stated explicitly in {\bf\cite{t}}, with a proof based on
{\bf\cite{du}}.  See {\bf\cite{sr2}} for a 
proof based on the ${\mathfrak {sl}}(2,\C)$ method
{\bf\cite{p}}. Very closely related 
results are shown in {\bf\cite{d2,st1,mp2,bvp}}. The
existence of an orthogonal SJB satisfying (\ref{sv}) has several
applications: in {\bf\cite{sr3}} we showed that the commutant of the
$GL(n,\Fq)$-action on $V(\Bq(n))$ block diagonalizes with respect to the
orthonormal basis given by the normalization of $J_q(n)$ and we used
(\ref{sv}) to make this block diagonalization explicit, thereby obtaining
a $q$-analog of the formula from {\bf\cite{s}} for explicit block
diagonalization of the commutant of the symmetric group action on $V(B(n))$.
This includes, as a special case, a formula for the eigenvalues of the
elements of the Bose-Mesner algebra
of the Grassmann scheme {\bf\cite{d1,cst}}. For other approaches to 
explicit block diagonalization see {\bf\cite{mp2,bvp}}, the latter of which
also gives applications to bounds on 
projective codes using semidefinite programming. In {\bf\cite{sr2}} we used
(\ref{sv}) to give a positive combinatorial formula for the number of
spanning trees of the $q$-analog of the $n$-cube and to show that
the Laplacian eigenvalues of the Grassmann graphs, known in principle since
{\bf\cite{d1}}, admit an elegant closed form. For another approach to
the Laplacian eigenvalues of the Grassmann graphs see {\bf\cite{mp2}}. At
the end of this paper we pose a bijective proof problem on spanning trees of
the Grassmann graphs.

From the point of view of applications (and especially that of polynomial
time computation) explicit construction of the basis
$J_q(n)$ is not important as even to write down $J_q(n)$ takes exponential
time. But it is of interest from a mathematical standpoint yielding useful
additional insight into the linear structure of the subspace lattice.
The situation is similar to bijective versus nonbijective proofs in
enumeration. Substituting $q=1$ in Theorem \ref{mt3} we recover the explicit orthogonal
SJB of $V(B(n))$ constructed in {\bf\cite{sr1}}. This basis was given a
representation theoretic characterization in {\bf\cite{sr1}}, namely, that
it is the canonically defined symmetric Gelfand-Tsetlin basis of $V(B(n))$.
Similarly, the basis $J_q(n)$ should also be studied from a representation 
theoretic viewpoint. We hope to return to this later.

{{\bf  \section {  Goldman-Rota recurrence}}}  

In this section we prove Theorem \ref{grv}. As stated in the introduction,
we identify $\Fq^k$, for $k<n$, with the subspace of $\Fq^n$ consisting  of
all vectors with the last $n-k$ components zero. We denote by $e_1,\ldots
,e_n$ the standard basis vectors of $\Fq^n$. So $\Bq(k)$ consists of all
subspaces of $\Fq^n$ contained in the subspace spanned by $e_1,\ldots ,e_k$.  

Define $\Aq(n)$ to be the
collection of all subspaces in $\Bq(n)$ not contained in the hyperplane
$\Fq^{n-1}$
$$\Aq(n) = \Bq(n) - \Bq(n-1) = \{ X\in \Bq(n) : X\not\subseteq \Fq^{n-1} \},\;n\geq 1.$$ 
For $1\leq k \leq n$, let $\Aq(n)_k$ denote the set of all subspaces in
$\Aq(n)$ with dimension $k$. We consider $\Aq(n)$ as an induced subposet of $\Bq(n)$. 

Define a map
$$\cH(n) : \Aq(n) \rar \Bq(n-1)$$
by $\cH(n)(X) = X\cap \Fq^{n-1}$, for $X\in \Aq(n)$. Define an equivalence
relation $\sim$ on $\Aq(n)$ by $X\sim Y$ iff $\cH(n)(X)=\cH(n)(Y)$. Denote the
equivalence class of $X\in \Aq(n)$ by $[X]$.

For $X\in \Bq(n-1)$, define $\wh{X}$ to be the subspace in $\Aq(n)$
spanned by $X$ and $e_n$.

\bl \label{el}
Let $X,Y \in \Aq(n)$ and $Z,T\in \Bq(n-1)$. Then

\noi (i) $\mbox{dim}\,\cH(n)(X)=\mbox{dim}\,X - 1$ and $\wh{\cH(n)(X)}\in
[X]$. 

\noi (ii) $Z\leq T$ iff $\wh{Z} \leq \wh{T}$.

\noi (iii) $Y$ covers $X$ iff $\cH(n)(Y)$ covers $\cH(n)(X)$ and 
$Y = \mbox{ span}\,(\cH(n)(Y)\cup \{v\})$ for any
$v\in X - \Fq^{n-1}$. 

\noi (iv) $| \cH(n)^{-1}(Z)| = q^l$, where $l= n-\mbox{dim}\,Z - 1$.
Thus, $|[X]|=q^{n-k}$, where $k=\mbox{ dim}\,X$.

\el
\pf (i), (ii), and (iii) are clear.

\noi (iv) Let $\cH(n)^{-1}(Z)=\{Y_1,\ldots ,Y_t\}$. Then $Y_i\cap Y_j = Z$,
$1\leq i\not= j \leq t$. Let $\mbox{dim}\,(Z)=m$. Now $|Y_i - Z| = q^{m+1} -
q^m$ for all $i$ and thus $t=\frac{q^n - q^{n-1}}{q^{m+1} - q^m}=q^{n-m-1}$.
$\Box$

We have an orthogonal decomposition
\beq \label{bod}
V(\Bq(n+1)) = V(\Bq(n)) \oplus V(\Aq(n+1)).
\eeq
We shall now give a canonical orthogonal decomposition of $V(\Aq(n+1))$.

Let $H(n+1, q)$ denote the subgroup of $GL(n+1,q)=GL(n+1,\Fq)$ consisting of all
matrices of the form
$$ \left[ \ba{cc}
           I & \ba{c} a_1 \\  \cdot \\ \cdot \\  a_n \ea \\
           0 \cdots 0 & 1
          \ea 
   \right],
$$
where $I$ is the $n\times n$ identity matrix.

The additive abelian group $\Fq^n$ is isomorphic to $H(n+1,q)$ via
$\phi : \Fq^n \rar H(n+1,q)$ given by
$$\phi \left(\left[ \ba{c} a_1 \\ \cdot \\ \cdot \\a_n \ea \right]\right)
\rar
\left[ \ba{cc}
           I & \ba{c} a_1 \\  \cdot \\ \cdot \\  a_n \ea \\
           0 \cdots 0 & 1
          \ea 
   \right].
$$   

There is a natural (left) action
of $H(n+1, q)$ on $\Aq(n+1)$ and $\Aq(n+1)_k$. For $X\in \Aq(n+1)$, let
$G_X \subseteq H(n+1,q)$ denote the stabilizer of $X$.

\bl \label{orl}
Let $X,Y\in \Aq(n+1)$. Then

\noi (i) The orbit of $X$ under the action of $H(n+1,q)$ is $[X]$.

\noi (ii) Suppose $Y$ covers $X$. Then the bipartite graph of the covering
relations between $[Y]$ and $[X]$ is regular with degrees $q$ (on the
$[Y]$ side) and $1$ (on the $[X]$ side). 
%Thus, for $W\in \Bq(n)$, we have
%\beq \label{succ}
%\sum_{S\sim \wh{W}} U_{n+1}(S) &=& q \left\{ \sum_{Z}
%\sum_{T\sim \wh{Z}} T \right\},
%\eeq 
%where the outer sum on the right hand side is over all $Z$ covering $W$ in
%$\Bq(n)$.

\noi (iii) Suppose $X\subseteq Y$. Then $G_X \subseteq G_Y$.
\el
\pf (i) This is clear.

\noi (ii) Since the action of $H(n+1,q)$ on $\Aq(n+1)$ is clearly order
preserving it follows that the bipartite graph in the statement is regular.
Let $Y'\in[Y]$ also cover $X$. Then $\cH(n+1)(Y')=\cH(n+1)(Y)$ and it follows
from Lemma \ref{el}(iii) that $Y=Y'$. So
the degree on the $[X]$ side is $1$. 
Let $\mbox{dim}\,(X)=k$. Then, by
Lemma \ref{el}(iv), $|[Y]|= q^{n-k}$ and $|[X]|=q^{n+1-k}$ and hence, by
regularity, the degree on the $[Y]$ side is $q$.

\noi (iii) We may assume that $Y$ covers $X$. If $X\not= X'\in [X]$ then
clearly $X\cap X' = \cH(n+1)(X)$. So, by part (ii) and Lemma
\ref{el}(iii), we can write $Y$ as a union
$$Y=X_1\cup X_2\cup \cdots X_q$$
of subspaces $X=X_1,\ldots , X_q \in [X]$ with $X_i \cap X_j =
\cH(n+1)(X),\;1\leq i\not= j\leq q$.

Now the stabilizer of all the
elements $X_1,\ldots ,X_q$ is $G_X$ (since $H(n+1,q)$ is commutative). It
follows that $G_X\subseteq G_Y$. $\Box$

Let $\cIq(n)$  
denote the set of all distinct irreducible characters (all of
degree $1$) of $H(n+1,q)$
and let $\cNq(n)$  
denote the set of all distinct nontrivial irreducible characters 
of $H(n+1,q)$.
 
Let $\psi_k$ (respectively, $\psi$) denote the character of the
permutation
representation of $H(n+1,q)$ on $V(\Aq(n+1)_k)$ (respectively,
$V(\Aq(n+1))$) corresponding to the left action. Clearly $\psi = \sum_{k=1}^{n+1} \psi_k$. 
Below $[,]$
denotes character inner product and the $q$-binomial coefficient $\qb{n}{k}$ is taken to
be zero when $n$ or $k$ is $< 0$.

\bt \label{mt1}
\noi (i) For $1\leq k \leq n+1$ and $g\in H(n+1,q)$ we have
$$\psi_k(g) = \left\{ \ba{ll} 
                   q^{n-k+1}\, \qb{n}{k-1} & \mbox{if }g=I,\\
                                   & \\
                                  q^{n-k+1}\, \qb{n-1}{k-2} & \mbox{if }g\not=I.
                                      \ea \right. $$
\noi (ii) Let $\chi \in \cIq(n)$ be the trivial character. Then
$[\chi , \psi_k ]= \qb{n}{k-1},\;1\leq k \leq n+1.$

\noi (iii) Let $\chi \in \cNq(n)$. Then
$[\chi , \psi_k ]= \qb{n-1}{k-1},\;1\leq k \leq n+1.$

\noi (iv) \beq \label{nsw}
\qb{n+1}{k} &=& \qb{n}{k} + \qb{n}{k-1} + (q^n - 1)\qb{n-1}{k-1},\;n,k\geq
1,\eeq
with $\qb{0}{k}=\delta(0,k)$ (Kronecker delta) and $\qb{n}{0}=1$. 
Note that (\ref{gri}) follows by summing over $k$.
\et
\pf (i) From Lemma \ref{el}(iv),
$$|\Aq(n+1)_k|=q^{n-k+1}\,\qb{n}{k-1}=\psi_k(I).$$ 
Now assume $g\in H(n+1,q),\,g\not= I$ and let $X\in \Aq(n+1)_k$. Let the
last column of $g$ be $(a_1,\ldots ,a_n,1)^t$ ($t$=transpose), where not all
the $a_i$'s are $0$. Now note that

\noi (a) For $b_1,\ldots , b_n \in \Fq$,
$$g(b_1,\ldots ,b_n,1)^t = (a_1+b_1,\ldots ,a_n+b_n,1)^t.$$

\noi (b) From item (a) above it follows that $gX=X$ iff the nonzero vector
$(a_1,\ldots ,a_n,0)^t \in \cH(n+1)(X)$.

\noi (c) From item (b) above, if $gX=X$ then $gY=Y$ for all $Y\in [X]$. Thus
$g$ either fixes all elements of $[X]$ or no elements.

\noi (d) The number of subspaces in $\Bq(n)_{k-1}$ containing the nonzero
vector $(a_1,\ldots ,a_n,0)^t$ is $\qb{n-1}{k-2}$.

It follows from items (b), (c), (d) above that $\psi_k(g) =
q^{n-k+1}\,\qb{n-1}{k-2}$.

\noi (ii) This follows from the well known result that the multiplicity of
the trivial representation in a permutation representation is the number of
orbits, which in the present case is $\qb{n}{k-1}$.

\noi (iii) Since $\chi$ is nontrivial we have $\ds{\sum_{g\in H(n+1,q)}
\ol{\chi(g)}} = 0$ and thus $\ds{\sum_{g\in H(n+1,q),\,g\not= I} 
\ol{\chi(g)} = -1}$. 

Thus (below the sum is over all $g\not= I$ in $H(n+1,q)$)
\beqn
[\chi,\psi_k] 
&=&\frac{1}{q^n}\left\{ \psi_k(I) + \sum_{g\not= I}
\ol{\chi(g)} \psi_k(g)\right\}\\
&=&\frac{1}{q^n}\left\{q^{n-k+1} \qb{n}{k-1} 
- q^{n-k+1} \qb{n-1}{k-2}  \right\}\\
&=& q^{-(k-1)}\left\{ \qb{n}{k-1} - \qb{n-1}{k-2} \right\}\\
&=& \qb{n-1}{k-1},
\eeqn
where in the last step we have used $q$-Pascal's triangle (see Section 1.7
in {\bf\cite{st2}})
$$\qb{n}{k-1} = \qb{n-1}{k-2} + q^{k-1} \qb{n-1}{k-1}.$$

\noi (iv) Let $1\leq k \leq n+1$. Restricting (\ref{bod}) to 
dimension $k$ we get the following
orthogonal decomposition
\beq 
V(\Bq(n+1)_k) = V(\Bq(n)_k) \oplus V(\Aq(n+1)_k).
\eeq
Splitting $V(\Aq(n+1)_k)$ into $H(n+1,q)$-irredicibles and taking dimensions
using parts (ii) and (iii) we get the result. The initial conditions are
easily verified. $\Box$

Let $W(0)$ (respectively, $W(0)_k$) denote the isotypical component of
$V(\Aq(n+1))$ (respectively, $V(\Aq(n+1)_k)$) corresponding to the trivial
representation of $H(n+1,q)$ and, for $\chi\in \cNq(n)$,
let $W(\chi)$ (respectively, $W(\chi)_k$) denote the isotypical component of
$V(\Aq(n+1))$ (respectively, $V(\Aq(n+1)_k)$) corresponding to the 
irreducible representation of $H(n+1,q)$ with character $\chi$.
We have the following orthogonal decompositions, the last of which is
canonical (note that $W(\chi)_{n+1}$ is the zero module, by Theorem
\ref{mt1}(iii)). 
\beq
W(0) &=& W(0)_1 \oplus \cdots \oplus W(0)_{n+1}, \\
W(\chi) &=& W(\chi)_1 \oplus \cdots \oplus W(\chi)_{n},\;\;\;\chi \in \cNq(n), \\
V(\Aq(n+1)) &=& W(0) \oplus \left(\oplus_{\chi \in \cNq(n)} W(\chi)\right).
\eeq
Since $U_{n+1}$ is $GL(n+1,q)$-linear, each of $W(0)$ and $W(\chi)$,
$\chi\in \cNq(n)$ is $U_{n+1}$-closed.

For $\chi \in \cIq(n)$, define the following element of the group algebra of
$H(n+1,q)$:
$$p(\chi)=\sum_g \ol{\chi(g)}\,g,$$
where the sum is over all $g\in H(n+1,q)$. For $1\leq k \leq n+1$, the map
\beq \label{proj} 
&p(\chi) : V(\Aq(n+1)_k) \rar V(\Aq(n+1)_k),& 
\eeq
given by $v\mapsto \sum_{g\in H(n+1,q)} \ol{\chi(g)}\,gv$,
is a nonzero multiple of the $H(n+1,q)$-linear projection onto
$W(\chi)_k$. Similarly for $p(\chi) : V(\Aq(n+1)) \rar V(\Aq(n+1))$.

\bl \label{cl} Let $X\in\Aq(n+1)$ and $\chi\in \cIq(n)$. Then $p(\chi)(X)\not=
0$ iff $\chi : G_X \rar \C^*$ is the trivial character of $G_X$.
\el
\pf Let $\{h_0=1,h_1,\ldots ,h_t\}$ be a set of distinct coset
representatives of $G_X$, i.e.,
$$H(n+1,q) = G_Xh_0 \cup G_Xh_1 \cup \cdots \cup G_Xh_t\;\;\mbox{
(disjoint union)}.$$
Write $[X]=\{X=X_0,X_1,\ldots ,X_t\}$ and assume without loss of generality
that $h_iX=X_i,\;0\leq i \leq t$. Note that $G_X$ is the stabilizer of all
the elements of $[X]$.

We have
\beq \nonumber
p(\chi)(X) &=& \sum_{g\in H(n+1,q)} \ol{\chi(g)}\,gX \\ \label{ci}
 &=& \left(\sum_{g\in G_X} \ol{\chi(g)}\right)\, X + 
               \sum_{i=1}^t \ol{\chi(h_i)}
\left(\sum_{g\in G_X} \ol{\chi(g)}\right)\, X_i .
\eeq
The result follows since $\sum_{g\in G_X} \ol{\chi(g)} = 0$ for every
nontrivial character of $G_X$. $\Box$

\bt \label{f} 

(i) Let $\chi\in \cIq(n),\;X,Y\in \Aq(n+1)$ with $X\sim Y$. Then
$p(\chi)(X)$ is a nonzero multiple of $p(\chi)(Y)$.

\noi (ii) Let $\chi\in \cIq(n)$. Then 
$\{ p(\chi)(\wh{X}) : X\in \Bq(n)_{k-1} \mbox{ with } p(\chi)(\wh{X})
\not= 0\}$ is a basis of $W(\chi)_k$, $1\leq k \leq n+1$.

\noi (iii) Let $\chi\in \cIq(n)$ and let $X,Y\in \Bq(n)$ with $X$ covering $Y$. Then
$$p(\chi)(\wh{X})\not= 0 \mbox{ implies } p(\chi)(\wh{Y})\not= 0.$$

\noi (iv) Define  $\theta_n : V(\Bq(n)) \rar W(0)$
by 
$$X\mapsto \sum_{Y\sim \wh{X}} Y,\;\,X\in \Bq(n).$$
Then $\theta_n$ is an isomorphism
of pairs $(V(\Bq(n)), qU_n)$ and $(W(0), U_{n+1})$ and
\beq \label{ul}
U_{n+1}(v) &=& U_n(v) + \theta_n(v),\,\,v\in V(\Bq(n)),\\
\label{ul1}
\langle \theta_n (w), \theta_n (v)\rangle &=& q^{n-k} \langle w,
v\rangle,\,\,w,v\in V(\Bq(n)_k),\;0\leq k \leq n. 
\eeq

\noi (v) Let $\chi\in \cNq(n)$.
From Theorem \ref{mt1} (iii) we have $\mbox{dim}\,W(\chi)_n = 1$. It thus
follows from part (ii) that there is a unique element 
$X=X(\chi)\in \Bq(n)_{n-1}$ such that
$p(\chi)(\wh{X})\not= 0$. Define
$ \lambda(\chi) : V(\Bq(X)) \rar W(\chi)$
by $$Y \mapsto p(\chi)(\wh{Y}),\;Y\in\Bq(X).$$ 
Then $\lambda(\chi)$ is an isomorphism of pairs
$(V(\Bq(X)), U_X)$ and $(W(\chi),U_{n+1})$ and satisfies
\beq 
\label{ul2}
\langle \lambda(\chi)(w), \lambda(\chi)(v)\rangle &=& q^{n+k} \langle w,
v\rangle,\,\,w,v\in V(\Bq(X)_k),\;0\leq k \leq n-1. 
\eeq

\noi (vi) For $X\in \Bq(n)_{n-1}$ the number of $\chi\in \cNq(n)$ such that
$p(\chi)(\wh{X})\not= 0$ is $q-1$.

\et
\pf (i) By Lemma \ref{orl}(i), $X=hY$ for some $h\in H(n+1,q)$. 
Then we have (below the sum is over all $g\in H(n+1,q)$)
$$
p(\chi)(X) = \sum_g \ol{\chi(g)}\,ghY 
           =  \sum_g \ol{\chi(gh^{-1})}\,gY
           =  \ol{\chi(h^{-1})}\,p(\chi)(Y).
$$
Since $\chi(h^{-1})\not= 0$ ($\chi$ being of degree 1) the result follows. 

\noi (ii) The map (\ref{proj}) is a projection onto $W(\chi)_k$, so 
${\cal G} = \{ p(\chi)(X) : X\in \Aq(n+1)_k\}$ spans $W(\chi)_k$. 
By part (i), the subset 
${\cal G'} = \{ p(\chi)(\wh{X}) : X\in \Bq(n)_{k-1} \mbox{ with }
p(\chi)(\wh{X})\not= 0\}$ also spans $W(\chi)_k$. 
Now, for distinct $X,Y\in \Bq(n)_{k-1}$, $p(\chi)(\wh{X})$ and
$p(\chi)(\wh{Y})$ have disjoint supports, so ${\cal G'}$ is a basis.

\noi (iii) This follows from Lemma \ref{orl}(iii) and Lemma \ref{cl}.

\noi (iv) 
By Theorem \ref{mt1} (ii) the dimensions of $V(\Bq(n))$ and
$W(0)$ are the same. For $X_1\not= X_2\in \Bq(n)$ the supports of
$\theta_n(X_1)$ and $\theta_n(X_2)$ are disjoint. It follows that
$\theta_n$ is a vector space isomorphism.  

Let $X\in\Bq(n)$ with $\mbox{dim}\,(X)=k$. Then (\ref{ul}) is
clear and $|[\wh{X}]|=q^{n-k}$
by Lemma \ref{el}(iv), showing (\ref{ul1}). 

We have (below the sum is over all $Z$ covering $X$ in $\Bq(n)$) 
\beqn
\theta_n(qU_n(X))&=&\theta_n\left(q\left(\sum_Z
Z\right)\right) \\
                        &=& q\,\left\{\sum_Z \sum_{Y\sim
\wh{Z}} Y\right\}.
\eeqn
Similarly (in the second step below $T$ varies over all subspaces covering
$Y$ and in the third step $Z$ varies over all subspaces in $\Bq(n)$ covering
$X$. We have used Lemma \ref{el}(ii) and Lemma \ref{orl}(ii) 
to go from the second to the third step)
\beqn
U_{n+1}(\theta_n(X))&=& 
U_{n+1}\left( \sum_{Y\sim \wh{X}}
Y\right)\\
&=&\sum_{Y\sim \wh{X}} \sum_T T\\
&=&q\,\left\{\sum_Z \sum_{Y\sim \wh{Z}} Y\right\}.
\eeqn

\noi (v) By part (iii) it follows that $\lambda(\chi)(Y)\not= 0$ for all
$Y\in \Bq(X)$. 
By Theorem \ref{mt1} (iii) the dimensions of $V(\Bq(X))$ and
$W(\chi)$ are the same. For $Y_1\not= Y_2\in \Bq(X)$ the supports of
$\lambda(\chi)(Y_1)$ and $\lambda(\chi)(Y_2)$ are disjoint. It follows that
$\lambda(\chi)$ is a vector space isomorphism.  

Now, for $Y\in\Bq(X)$, we have (below the sum is over all $Z$ covering $Y$
in $\Bq(X)$)
\beqn
\lambda(\chi)(U_X(Y)) &=& 
\lambda(\chi)\;\left(\sum_Z Z\right)\\
             &=& 
\sum_Z p(\chi)(\wh{Z}). 
\eeqn

Let $Y\in \Bq(X)$. Before calculating $U_{n+1}\lambda(\chi)(Y)$ we make the
following observation. By Lemma \ref{el}(ii) every element covering $\wh{Y}$
is of the form $\wh{Z}$, for some $Z$ covering $Y$ in $\Bq(n)$. Suppose
$Z\in \Bq(n) - \Bq(X)$. Since $\mbox{dim}\,(W(\chi))=G_q(n-1)$ (by Theorem
\ref{mt1}(iii)), it follows by parts (ii) and (iii) that
$p(\chi)(\wh{Z})=0$.

We now calculate $U_{n+1}\lambda(\chi)(Y)$. In the second step below we have
used the fact that $U_{n+1}$ is $H(n+1,q)$-linear and in the third step,
using the observation in the paragraph above, we may restrict the sum to all
$Z$ covering $Y$ in $\Bq(X)$.

We have
\beqn
U_{n+1}(\lambda(\chi)(Y))&=& U_{n+1}\left( p(\chi)(\wh{Y}) \right)\\
    &=& p(\chi)(U_{n+1}(\wh{Y}))\\ 
    &=& \sum_Z p(\chi)(\wh{Z}). 
\eeqn
We will now show that $\pr p(\chi)(\wh{Y})\pr = \sqrt{q^{n+k}}$ if $Y\in
\Bq(X)$ with $\mbox{dim}\,(Y)=k$. This will prove (\ref{ul2}).

By Lemma \ref{el}(iv), $|[\wh{Y}]|=q^{n+1-(k+1)}=q^{n-k}$ and so
$|G_{\wh{Y}}| = q^k$. It now follows from (\ref{ci}) (since the restriction
of $\chi$ to $G_{\wh{Y}}$ must be trivial) that 
$p(\chi)(\wh{Y})=\sqrt{q^{n-k}(q^k)^2}=\sqrt{q^{n+k}}.$ 

\noi (vi) By Theorem \ref{mt1}(iii), $\sum_{\chi \in \cNq(n)}
\mbox{dim}\,(W(\chi)_n)=q^n - 1=(q-1)(1+q+\cdots +q^{n-1})$. 
Now $|\Bq(n)_{n-1}|=\qb{n}{n-1}=1+q+\cdots +q^{n-1}$ and by Lemma
\ref{el}(iv), $|[\wh{X}]|=q$ for $X\in \Bq(n)_{n-1}$. Since, for $\chi\in \cNq(n)$
and $X\in\Bq(n)_{n-1}$, the support of $p(\chi)(\wh{X})$ is contained in
$[\wh{X}]$ and $p(\chi)(\wh{X})$ is orthogonal to $p(\pi)(\wh{X})$ (where $\pi$
is the trivial character), the result now follows by part (ii). $\Box$

To use Theorem \ref{f} for computations we need the character table of
$H(n+1,q)$, which is easy to write down explicitly since $H(n+1,q)$ is
direct sum of $n$ cyclic groups of order $q$. 
We now give a small example to illustrate part (v) of Theorem
\ref{f}.

\noi {\bf Example} Let $q=3$, $n=2$, and $\omega = e^{2\pi \mbox{i}/3}$. 
Consider $A_3(3)$ with the
$H(3,3)$-action. We write the elements of ${\mathbb F}_3$ as $\{0,1,2\}$
and define $\chi \in {\cal N}_3(2)$ by
$\chi(\phi((a_1,a_2)) = \omega^{a_1 + 2a_2}$, where $a_1+2a_2$ is computed
as an integer.

We have
\beqn
p(\chi) &=&
  \left[ \begin{array}{ccc} 1&0&0\\0&1&0\\0&0&1 \end{array} \right] 
  + \omega^2 \left[ \begin{array}{ccc} 1&0&0\\0&1&1\\0&0&1 \end{array} \right] 
  + \omega \left[ \begin{array}{ccc} 1&0&0\\0&1&2\\0&0&1 \end{array}
\right]\\ &&\\ 
 && + \;\omega \left[ \begin{array}{ccc} 1&0&1\\0&1&0\\0&0&1 \end{array} \right] 
  + \left[ \begin{array}{ccc} 1&0&1\\0&1&1\\0&0&1 \end{array} \right] 
  + \omega^2 \left[ \begin{array}{ccc} 1&0&1\\0&1&2\\0&0&1 \end{array}
\right]\\&&\\
 && +\; \omega^2 \left[ \begin{array}{ccc} 1&0&2\\0&1&0\\0&0&1 \end{array} \right] 
  + \omega \left[ \begin{array}{ccc} 1&0&2\\0&1&1\\0&0&1 \end{array} \right] 
  + \left[ \begin{array}{ccc} 1&0&2\\0&1&2\\0&0&1 \end{array}
\right]
\eeqn
Given a finite set of vectors $v_1,v_2,\ldots ,v_m$ we shall denote the
subspace spanned by them by $\cS(v_1,\ldots ,v_m)$.

The four subspaces in $B_3(2)_1$ are
$$
X_1=\cS\left(\left[ \ba{c} 1\\0\\0 \ea \right]\right),\;
X_2=\cS\left(\left[ \ba{c} 0\\1\\0 \ea \right]\right),\;
X_3=\cS\left(\left[ \ba{c} 1\\1\\0 \ea \right]\right),\;
X_4=\cS\left(\left[ \ba{c} 2\\1\\0 \ea \right]\right).
$$
It can be checked that
\beqn
p(\chi)(\wh{X_1}) &=&
\cS\left(\left[\ba{c}1\\0\\0\ea\right],\;\left[\ba{c}0\\0\\1\ea\right]\right)
+ \omega^2 
\cS\left(\left[\ba{c}1\\0\\0\ea\right],\;\left[\ba{c}0\\1\\1\ea\right]\right)
+ \omega
\cS\left(\left[\ba{c}1\\0\\0\ea\right],\;\left[\ba{c}0\\2\\1\ea\right]\right)\\
&&\\
&&
+ \;\omega 
\cS\left(\left[\ba{c}1\\0\\0\ea\right],\;\left[\ba{c}1\\0\\1\ea\right]\right)
+ 
\cS\left(\left[\ba{c}1\\0\\0\ea\right],\;\left[\ba{c}1\\1\\1\ea\right]\right)
+ \omega^2
\cS\left(\left[\ba{c}1\\0\\0\ea\right],\;\left[\ba{c}1\\2\\1\ea\right]\right)\\
&&\\
&&
+ \;\omega^2 
\cS\left(\left[\ba{c}1\\0\\0\ea\right],\;\left[\ba{c}2\\0\\1\ea\right]\right)
+ \omega
\cS\left(\left[\ba{c}1\\0\\0\ea\right],\;\left[\ba{c}2\\1\\1\ea\right]\right)
+ 
\cS\left(\left[\ba{c}1\\0\\0\ea\right],\;\left[\ba{c}2\\2\\1\ea\right]\right)
\eeqn
Using $\omega^3=1$ and $1+\omega+\omega^2=0$ we see that
$p(\chi(\wh{X_1}))=0$.

Similarly we can check that 
\beqn
p(\chi)(\wh{X_2}) &=&
\cS\left(\left[\ba{c}0\\1\\0\ea\right],\;\left[\ba{c}0\\0\\1\ea\right]\right)
+ \omega^2 
\cS\left(\left[\ba{c}0\\1\\0\ea\right],\;\left[\ba{c}0\\1\\1\ea\right]\right)
+ \omega
\cS\left(\left[\ba{c}0\\1\\0\ea\right],\;\left[\ba{c}0\\2\\1\ea\right]\right)\\
&&\\
&&
+\; \omega 
\cS\left(\left[\ba{c}0\\1\\0\ea\right],\;\left[\ba{c}1\\0\\1\ea\right]\right)
+ 
\cS\left(\left[\ba{c}0\\1\\0\ea\right],\;\left[\ba{c}1\\1\\1\ea\right]\right)
+ \omega^2
\cS\left(\left[\ba{c}0\\1\\0\ea\right],\;\left[\ba{c}1\\2\\1\ea\right]\right)\\
&&\\
&&
+ \;\omega^2 
\cS\left(\left[\ba{c}0\\1\\0\ea\right],\;\left[\ba{c}2\\0\\1\ea\right]\right)
+ \omega
\cS\left(\left[\ba{c}0\\1\\0\ea\right],\;\left[\ba{c}2\\1\\1\ea\right]\right)
+ 
\cS\left(\left[\ba{c}0\\1\\0\ea\right],\;\left[\ba{c}2\\2\\1\ea\right]\right)
\eeqn
and
\beqn
p(\chi)(\wh{X_4}) &=&
\cS\left(\left[\ba{c}2\\1\\0\ea\right],\;\left[\ba{c}0\\0\\1\ea\right]\right)
+ \omega^2 
\cS\left(\left[\ba{c}2\\1\\0\ea\right],\;\left[\ba{c}0\\1\\1\ea\right]\right)
+ \omega
\cS\left(\left[\ba{c}2\\1\\0\ea\right],\;\left[\ba{c}0\\2\\1\ea\right]\right)\\
&&\\
&&
+ \;\omega 
\cS\left(\left[\ba{c}2\\1\\0\ea\right],\;\left[\ba{c}1\\0\\1\ea\right]\right)
+ 
\cS\left(\left[\ba{c}2\\1\\0\ea\right],\;\left[\ba{c}1\\1\\1\ea\right]\right)
+ \omega^2
\cS\left(\left[\ba{c}2\\1\\0\ea\right],\;\left[\ba{c}1\\2\\1\ea\right]\right)\\
&&\\
&&
+ \;\omega^2 
\cS\left(\left[\ba{c}2\\1\\0\ea\right],\;\left[\ba{c}2\\0\\1\ea\right]\right)
+ \omega
\cS\left(\left[\ba{c}2\\1\\0\ea\right],\;\left[\ba{c}2\\1\\1\ea\right]\right)
+ 
\cS\left(\left[\ba{c}2\\1\\0\ea\right],\;\left[\ba{c}2\\2\\1\ea\right]\right)
\eeqn
are both equal to $0$ and
\beqn
p(\chi)(\wh{X_3}) &=&
\cS\left(\left[\ba{c}1\\1\\0\ea\right],\;\left[\ba{c}0\\0\\1\ea\right]\right)
+ \omega^2 
\cS\left(\left[\ba{c}1\\1\\0\ea\right],\;\left[\ba{c}0\\1\\1\ea\right]\right)
+ \omega
\cS\left(\left[\ba{c}1\\1\\0\ea\right],\;\left[\ba{c}0\\2\\1\ea\right]\right)\\
&&\\
&&
+ \;\omega 
\cS\left(\left[\ba{c}1\\1\\0\ea\right],\;\left[\ba{c}1\\0\\1\ea\right]\right)
+ 
\cS\left(\left[\ba{c}1\\1\\0\ea\right],\;\left[\ba{c}1\\1\\1\ea\right]\right)
+ \omega^2
\cS\left(\left[\ba{c}1\\1\\0\ea\right],\;\left[\ba{c}1\\2\\1\ea\right]\right)\\
&&\\
&&
+ \;\omega^2 
\cS\left(\left[\ba{c}1\\1\\0\ea\right],\;\left[\ba{c}2\\0\\1\ea\right]\right)
+ \omega
\cS\left(\left[\ba{c}1\\1\\0\ea\right],\;\left[\ba{c}2\\1\\1\ea\right]\right)
+ 
\cS\left(\left[\ba{c}1\\1\\0\ea\right],\;\left[\ba{c}2\\2\\1\ea\right]\right)
\eeqn
is $\not= 0$.

We have now proved most of Theorem \ref{grv} except for one small part. Let
$X\in \Bq(n)_{n-1}$. The pairs $(V(\Bq(X)),U_X)$ and $(V(\Bq(n-1)),U_{n-1})$
are clearly isomorphic with many possible isomorphisms. We now define a
canonical isomorphism, based on the concept of a matrix in Schubert normal
form.

A $n\times k$ matrix $M$ over $\Fq$ is in {\em Schubert normal form}
 (or, {\em column reduced echelon form}) provided

\noi (i) Every column is nonzero.

\noi (ii) The first nonzero entry in every column is a $1$. Let the first
nonzero entry in column $j$ occur in row $r_j$.

\noi (iii) We have $r_1 < r_2 < \cdots < r_k$ and the submatrix of $M$
formed by the rows $r_1,r_2,\ldots ,r_k$ is the $k\times k$ identity matrix.

It is well known 
that every $k$ dimensional subspace
of $\Fq^n$ is the column space of a unique $n\times k$ matrix in Schubert
normal form (see Proposition 1.7.3 in {\bf\cite{st2}} where the
discussion is in terms of the row space).

Let $X\in \Bq(n)_{n-1}$ and let $M(X)$ be the $n\times (n-1)$ matrix in
Schubert normal form with column space $X$. The map $\tau(X) : \Fq^{n-1}
\rar X$ given by $e_j \mapsto \mbox{ column $j$ of }X$ is clearly a linear
isomorphism and this isomorphism gives rise to an isomorphism
$$\mu(X) : V(\Bq(n-1)) \rar V(\Bq(X))$$
of pairs $(V(\Bq(n-1)),U_{n-1})$ and $(V(\Bq(X)),U_X)$ given by
$\mu(X)(Y)=\tau(X)(Y),\;Y\in\Bq(n-1)$.

\noi \pf (of Theorem \ref{grv}) It is convenient to write the orthogonal
decomposition (\ref{odsd}) as follows
\beq \label{od}
V(\Bq(n+1)) &=& V(\Bq(n)) \oplus 
W(0) \oplus \left(\oplus_{\chi \in \cNq(n)} W(\chi)\right).
\eeq
Note that $|\cNq(n)|=q^n - 1$.

\noi (i) We have already showed that each of $W(0)$ and $W(\chi),\;\chi\in\cNq(n)$ is
$U_{n+1}$-closed. The rank sets of $W(0)$ and $W(\chi)$ are also easily seen
to be as stated.

\noi (ii) This follows from Theorem \ref{f}(iv).

\noi (iii) Let $\chi\in\cNq(n)$. From Theorem \ref{f}(v) there is a unique
$X\in\Bq(n)_{n-1}$ with $p(\chi)(\wh{X})\not= 0$. It now follows, again by
Theorem \ref{f}(v), that $\gamma_{n-1}(\chi)=\lambda(\chi)\mu(X)$ is an isomorphism of
pairs $(V(\Bq(n-1)),U_{n-1})$ and $(W(\chi),U_{n+1})$ satisfying (\ref{ti2}). 
$\Box$

{{\bf  \section {  Orthogonal symmetric Jordan basis}}}  

In this section we prove Theorem \ref{mt3} and give an application to the
Grassmann scheme. We also pose a bijective proof problem on the Grassmann
graphs.

\pf (of Theorem \ref{mt3}) The proof is by induction on $n$, the result
being clear for $n=0,1$. 

Let $\chi\in \cNq(n)$ and let $(x_k,\ldots
,x_{n-1-k})$ be a SJC in $J_q(n-1)$ 
starting at rank $k$ and ending at rank $n-1-k$. Then, 
by Theorem \ref{grv}, applied to the decomposition (\ref{od}),
$(y_{k+1},\ldots ,y_{n-k})$, where $y_{u+1} = 
\gamma_{n-1}(\chi)(x_u),\,k\leq u \leq n-1-k$, is a SJC in $W(\chi)$ (with respect to
$U_{n+1}$) starting at rank $k+1$, 
ending at rank $n-k$. By the induction hypothesis we have,
for $k+1\leq u \leq n-k$, 
$$
\frac{\pr y_{u+1} \pr}{\pr y_u \pr} =
\frac{\sqrt{q^{n+u}}\pr x_{u} \pr}{\sqrt{q^{n+u-1}}\pr x_{u-1} \pr} =
%\sqrt{q}\,\sqrt{q^k\kq{u-k}\kq{n-k-u}} =
\sqrt{q^{k+1}\kq{u+1-(k+1)}\kq{n-k-u}}.
$$
Doing the above procedure for every SJC in $J_q(n-1)$ we get an orthogonal
SJB of $W(\chi)$ satisfying (\ref{sv}). Note that, by definition of
$\lambda(\chi)$, if the coefficients (in the standard basis) of the vectors
in $J_q(n-1)$ were integral multiples of $q\mbox{th}$ roots of unity then so
will be the coefficients of the vectors in the SJB of $W(\chi)$. 
Similarly, doing the above procedure for every $\chi\in\cNq(n)$ we get an 
orthogonal SJB, with respect to $U_{n+1}$, of $\oplus_{\chi\in \cNq(n)}
W(\chi)$ satisfying (\ref{sv}).

Now we consider the subspace $V(\Bq(n)) \oplus W_0$. Let $(x_k,\ldots
,x_{n-k})$ be a SJC in $J_q(n)$, starting at rank $k$ and ending at rank
$n-k$, and satisfying (\ref{sv}). 
Set $\ol{x_u}=\theta_n(x_u),\;k\leq
u \leq n-k$.  
Then, by Theorem
\ref{grv}, $(w_{k+1},\ldots ,w_{n-k+1})$, where $w_{u+1} =
q^{u-k}\,\ol{x_u},\,k\leq u \leq n-k$ is a graded Jordan chain
in $W_0$ (with respect to $U_{n+1}$), starting at rank $k+1$ and ending at
rank $n-k+1$. We have $U_{n+1}(q^{u-k}\,\ol{x_u})=q^{u+1-k}\,\ol{x_{u+1}}$ and
so
\beq \label{u1}
U_{n+1}(\ol{x_u})=q\,\ol{x_{u+1}},\;k\leq u < n-k.
\eeq
Also we have
\beq \label{bar}
\langle \ol{x_u} , \ol{x_u} \rangle = q^{n-u}\, \langle x_u , x_u \rangle
,\;k\leq u \leq n-k.
\eeq
For convenience we define $x_{k-1}=\ol{x_{k-1}}=x_{n+1-k}=0$. Note that
(\ref{bar}) also holds for $u=k-1$.

Now, by (\ref{ti}), we have, for $k\leq u \leq n-k$,
\beq \label{up}
&U_{n+1}(x_u) \,=\,  x_{u+1} + \theta_n(x_u)  
           \,=\, x_{u+1} + \ol{x_u}
\eeq 
Let $Z$ be the subspace spanned by $\{x_k ,\ldots ,x_{n-k}\}$ and
$\{\ol{x_k}, \ldots ,\ol{x_{n-k}}\}$. Clearly, by (\ref{u1}) and (\ref{up}), 
$Z$ is $U_{n+1}$-closed. We
shall now get an orthogonal SJB of $Z$ satisfying (\ref{sv}) by taking
linear combinations of the vectors $\{x_k,\ldots ,x_{n-k}\}$ and
$\{\ol{x_k}, \ldots ,\ol{x_{n-k}}\}$. 

We consider two cases:

\noi (a) $k=n-k$ : By (\ref{up}), $(x_k, \ol{x_k})$ is an
orthogonal SJB of $Z$ going from rank $k$ to rank $k+1$. We have, from
(\ref{bar}),
$$ \frac{\langle \ol{x_k} , \ol{x_k} \rangle}
{\langle x_k , x_k \rangle}= q^k,$$
and thus (\ref{sv}) is satisfied.

\noi (b) $k < n-k$ : Define the following vectors in $Z$.
\beqn
y_l &=& x_l + \kq{l-k}\,\;\ol{x_{l-1}},\;\;k\leq l \leq
n+1-k, \\
z_l &=& - q^n \,x_l + q^{l+k-1}\;\kq{n-l-k+1}\,\;
\ol{x_{l-1}},\;\;k+1\leq l \leq
n-k.
\eeqn
Note that, using the induction hypothesis, the coefficients of $y_l$, $z_l$
are also integral multiples of $q\mbox{th}$ roots of unity.
We claim that $(y_k,\ldots ,y_{n+1-k})$ and $(z_{k+1},\ldots ,z_{n-k})$ form
an orthogonal SJB of $Z$ satisfying (\ref{sv}).

We check orthogonality first, for which we need to show that 
$\langle y_l,z_l\rangle =0$
for $k+1\leq l \leq n-k$. 
Clearly $\langle x_l , \ol{x_{l-1}} \rangle =
0$ for $k+1\leq l \leq n-k$. Thus
$$ \langle y_l , z_l \rangle = - q^n \langle x_l , x_l \rangle
                               +
q^{k+l-1}\;\kq{l-k}\;\kq{n-l-k+1}\;\,\langle \ol{x_{l-1}} ,
\ol{x_{l-1}} \rangle .
$$
By the induction hypothesis $\langle x_l , x_l \rangle = q^k
\;\kq{l-k}\;\kq{n-l-k+1}\;\,\langle x_{l-1} , x_{l-1} \rangle$ and by
(\ref{bar}) $\langle \ol{x_{l-1}} , \ol{x_{l-1}} \rangle =
q^{n+1-l}\;\,\langle x_{l-1} , x_{l-1} \rangle$.
Thus $\langle y_l , z_l \rangle = 0$.

Now we check the Jordan chain condition. Using (\ref{u1}) and 
(\ref{up}), we have, for
$k\leq l < n+1-k$,
\beqn
U_{n+1}(y_l)&=& x_{l+1} +
\left(1+\,q\kq{l-k}\right)\;\ol{x_l} \\
&=& x_{l+1} + \kq{l+1-k}\;\ol{x_l} \\
&=& y_{l+1},
\eeqn
and clearly $U_{n+1}(y_{n+1-k})=0$.

Similarly, for $k+1\leq l \leq n-k$,
\beqn
U_{n+1}(z_l)&=& - q^n x_{l+1} +
q^{l+k}\;\left(\kq{n-l-k+1} - q^n\,\right)\;\ol{x_l} \\
&=& - q^n x_{l+1} + q^{l+k}\;\kq{n-l-k}\;\ol{x_l} \\
&=& z_{l+1}.
\eeqn
Note that $z_{n-k+1}=0$.

Now we check that condition (\ref{sv}) holds. 
For $k\leq u < n+1-k$ 
we have by the induction hypothesis (in the second step below we have used
(\ref{bar}). Note the second term in the denominator after the 
fourth step below. This is a fraction with a term $\kq{u-k}$ in the
denominator, which is zero for $u=k$. This is permissible here because of
the presence of the factor $\kq{u-k}^2$ in the numerator)
\beqn
\frac{\langle y_{u+1},y_{u+1}\rangle}{\langle y_u,y_u\rangle }&=&
\frac{\langle x_{u+1},x_{u+1}\rangle + \kq{u+1-k}^2 \;
\langle \ol{x_u} ,\ol{x_u} \rangle}
{\langle x_{u},x_{u}\rangle + \kq{u-k}^2 \;
\langle \ol{x_{u-1}} ,\ol{x_{u-1}} \rangle}\\
&=&\frac{\langle x_{u+1},x_{u+1}\rangle + q^{n-u}\;\kq{u+1-k}^2 \;
\langle x_u , x_u  \rangle}
{\langle x_{u},x_{u}\rangle + q^{n-u+1}\;\kq{u-k}^2 \;
\langle x_{u-1} ,x_{u-1} \rangle}\\
&=&\frac{ \frac{\langle x_{u+1},x_{u+1}\rangle}{\langle x_u , x_u \rangle }
           + q^{n-u}\;\kq{u+1-k}^2} 
{ 1 + q^{n-u+1}\;\kq{u-k}^2 \;
\frac{\langle x_{u-1} ,x_{u-1}\rangle }{\langle x_u , x_u \rangle}}\\
&=&\frac{ q^k \kq{u+1-k}\,\kq{n-k-u}
           + q^{n-u}\;\kq{u+1-k}^2} 
{ 1 + \frac{q^{n-u+1}\;\kq{u-k}^2}
           {q^k \kq{u-k}\;\kq{n-k-u+1}}}\\
&=& q^k\;\kq{u+1-k}\;\kq{n-k-u+1}\;\left(
 \frac{\kq{n-k-u} + q^{n-k-u}\;\kq{u+1-k}}{\kq{n-k-u+1} +
q^{n-k-u+1}\;\kq{u-k}}\right)\\
&=& q^k\;\kq{u+1-k}\;\kq{n-k-u+1}.
\eeqn
Similarly, for $k+1\leq u < n-k$, we have
\beqn
\frac{\langle z_{u+1},z_{u+1}\rangle}{\langle z_u,z_u\rangle }&=&
\frac{q^{2n}\langle x_{u+1},x_{u+1}\rangle + q^{2u+2k}\;\kq{n-u-k}^2 \;
\langle \ol{x_u} ,\ol{x_u} \rangle}
{q^{2n}\langle x_{u},x_{u}\rangle + q^{2u+2k-2}\;\kq{n-u-k+1}^2 \;
\langle \ol{x_{u-1}} ,\ol{x_{u-1}} \rangle}\\
&=&\frac{\langle x_{u+1},x_{u+1}\rangle + q^{2k-n+u}\;\kq{n-u-k}^2 \;
\langle x_u , x_u  \rangle}
{\langle x_{u},x_{u}\rangle + q^{2k-n+u-1}\;\kq{n-u-k+1}^2 \;
\langle x_{u-1} ,x_{u-1} \rangle}\\
&=&\frac{ \frac{\langle x_{u+1},x_{u+1}\rangle}{\langle x_u , x_u \rangle }
           + q^{2k-n+u}\;\kq{n-u-k}^2} 
{ 1 + q^{2k-n+u-1}\;\kq{n-u-k+1}^2 \;
\frac{\langle x_{u-1} ,x_{u-1}\rangle }{\langle x_u , x_u \rangle}}\\
&=&\frac{ q^k \kq{u+1-k}\,\kq{n-k-u}
           + q^{2k-n+u}\;\kq{n-u-k}^2} 
{ 1 + \frac{q^{2k-n+u-1}\;\kq{n-u-k+1}^2}
           {q^k \kq{u-k}\;\kq{n-k-u+1}}}\\
&=& q^{k+1}\;\kq{u-k}\;\kq{n-k-u}\;\left(
 \frac{\kq{n-k-u} + q^{n-k-u}\;\kq{u+1-k}}{\kq{n-k-u+1} +
q^{n-k-u+1}\;\kq{u-k}}\right)\\
&=& q^{k+1}\;\kq{u-k}\;\kq{n-k-u}.
\eeqn
Since $\theta_n$ is an isomorphism, 
doing the procedure above for every SJC in $J_q(n)$ we get an orthogonal SJB
of $V(\Bq(n)) \oplus W(0)$ satisfying (\ref{sv}).
That completes the proof. $\Box$

We now consider the application of Theorem \ref{mt3}
to the Bose-Mesner algebra of the Grassmann scheme of $m$-subspaces. 
For convenience we assume $0\leq m \leq n/2$. 
We do
not define this algebra here but instead work with the well known
characterization that it equals the commutant of the $GL(n,q)$-action on
$V(\Bq(n)_m)$. For the proof of the following result 
see Chapter 29 of {\bf\cite{jl}} where 
the $q=1$ case is proven. The same
proof works in general. 
\bt \label{iirep}
Let $0\leq m \leq n/2$. Then $V(\Bq(n)_m)$ is a multiplicity free
$GL(n,q)$-module with $m+1$ distinct irreducible summands. 
\et
Thus
$\mbox{End}_{\,GL(n,q)}(V(\Bq(n)_m)))$ is a commutative $*$-algebra with
dimension $m+1$ and so can be unitarily diagonalized. 

\bt Let $0\leq m \leq n/2$. Define
$$J_q(n,m) = \{ v\in J_q(n) : r(v)=m\}.$$
Then $J_q(n,m)$ is a common  orthogonal eigenbasis for the elements of
$\mbox{End}_{\,GL(n,q)}(V(\Bq(n)_m)))$.
\et
\pf 
For $i=0,1,\ldots ,m$ and  $k=0,1,\ldots ,i$ define
\beqn
J_q(n,i,k) &=& \{ v\in J_q(n) : r(v)=i \mbox{ and the Jordan chain
containing $v$}\\ 
\nonumber && \mbox{$\;\;\;\;\;\;\;\;\;\;\;\;\;\;\;\;\;\;\;\;\;$
                   starts at rank $k$}\}.
\eeqn 
Let $W_q(n,i,k)$ be the subspace spanned by $J_q(n,i,k)$. Then we have
an orthogonal direct sum decomposition
\beq \label{irred}
V(\Bq(n)_i) &=& \oplus_{k=0}^i W_q(n,i,k).
\eeq
Clearly $\mbox{dim}(W_q(n,i,k)) = \qb{n}{k} - \qb{n}{k-1}$. 

We shall now show that, for $i=0,1,\ldots ,m$,
$W_q(n,i,k),\;k=0,1,\ldots ,i$ are
$GL(n,q)$-submodules of $V(\Bq(n)_i)$. We do this by induction on $i$, the
case $i=0$ being clear.

Assume inductively that   
$W_q(n,i-1,0),\ldots ,W_q(n,i-1,i-1)$ are $GL(n,q)$-submodules, 
where $i < m$. 
Since $U_n$ is $GL(n,q)$-linear, 
$U_n(W_q(n,i-1,j))=W_q(n,i,j)$, $0\leq j \leq i-1$ are
$GL(n,q)$-submodules. Now consider $W_q(n,i,i)$.
Let $u\in W_q(n,i,i)$ and $g\in GL(n,q)$. Since $U_n$ is
$GL(n,q)$-linear we have $U_n^{n-2i+1}(g u)=g U_n^{n-2i+1}(u)=0$.
It follows that
$g u \in W_q(n,i,i)$.

We now have from Theorem \ref{iirep} that (\ref{irred}) is the
decomposition of $V(\Bq(n)_i)$ into distinct irreducible modules. The result
follows. $\Box$

Using (\ref{sv}) we can also determine the eigenvalues of the elements of 
$\mbox{End}_{\,GL(n,q)}(V(\Bq(n)_m)))$. More generally, we can explicitly
block diagonalize  
$\mbox{End}_{\,GL(n,q)}(V(\Bq(n)))$. We refer to {\bf\cite{sr3}} for
details.

Finally, we pose a bijective proof problem on the spanning trees of the
Grassmann and Johnson graphs. Actually, this application only requires the
existence of an orthogonal SJB satisfying (\ref{sv}) and not the actual
construction from the present paper.

The number of spanning trees of a graph $G$ is called the {\em
complexity} of $G$ and denoted $c(G)$. 
The number of {\em rooted spanning trees} (i.e., a spanning
tree plus a choice of a vertex as a root) of $G$ is denoted
$\oc{c}(G)$.  

Let $0\leq m \leq n/2$. The {\em Johnson graph} $C(n,m)$ is
defined to be the graph
with $B(n)_m$, the set of all subsets in $B(n)$ of cardinality
$m$, 
as the vertex set and with two vertices $X,Y\in B(n)_m$
connected by an edge iff $|X\cap Y| = m-1$.

Let
$0\leq m \leq n/2$. The {\em Grassmann graph} $\Cq(n,m)$ is
defined to be the graph
with vertex set $\Bq(n)_m$, and with two vertices $X,Y\in \Bq(n)_m$ 
connected by an edge iff $\mbox{dim}(X\cap Y) = m-1$. 

Let $T_q(n,m)$ and $T(n,m)$ denote, respectively, the set of rooted spanning
trees of
$\Cq(n,m)$ and $C(n,m)$.

For $X\in \Bq(n)_k, \; X' \in \Bq(n)_{k-1},\;1\leq k \leq n$ define
\beqn
{\cal U}{\cal D}(X) &=& \{ (Y,Z) \in \Bq(n)_{k-1} \times \Bq(n)_k \;|\; X
\supseteq Y \subseteq Z \},\\
{\cal D}{\cal U}(X') &=& \{ (Y',Z') \in \Bq(n)_{k} \times \Bq(n)_{k-1} \;|\;
X'
\subseteq Y' \supseteq Z' \}.\eeqn
\bt \label{gg}
Let $0\leq m \leq n/2$. The sets
$$T_q(n,m) \times \prod_{X\in \Bq(n)_{m-1}} {\cal D}{\cal U}(X)
\;\;\;\;\;\;\mbox{\em{and}}\;\;\;\;\;\; T_q(n,m-1) \times  
\prod_{X\in \Bq(n)_m} {\cal U}{\cal D}(X)$$
have the same cardinality. 
\et
\pf We give an algebraic proof.
For $X\in \Bq(n)_k, \; X' \in \Bq(n)_{k-1},\;1\leq k \leq n$ note that
$$ |{\cal U}{\cal D}(X)| = |{\cal D}{\cal U}(X')| = \kq{k}\kq{n-k+1}.$$
Now, using the existence of an orthogonal SJB of $V(\Bq(n))$ satisfying
(\ref{sv}) it was proved in {\bf\cite{sr2}} that the Laplacian eigenvalues
of $C_q(n,m)$ are $\kq{k}\kq{n-k+1},\;k=0,1,\ldots ,m$ with respective
multiplicities $\qb{n}{k}-\qb{n}{k-1}$. It now follows from the matrix-tree
theorem (see {\bf\cite{bh}}) that
\beqn 
\oc{c}(\Cq(n,m)) &=&
 {\ds{\prod_{k=1}^{m}}}
  \left( \kq{k}\kq{n-k+1} \right)^{\qb{n}{k}- \qb{n}{k-1}}.
\eeqn 
It follows that the sets in the statement of the theorem have the same
cardinality. $\Box$

The following result is an immediate corollary of the theorem above. We use
similar notations as above.
\bt \label{jg}
Let $0\leq m \leq n/2$. The sets
$$T(n,m) \times \prod_{X\in B(n)_{m-1}} {\cal D}{\cal U}(X)
\;\;\;\;\;\;\mbox{\em{and}}\;\;\;\;\;\; T(n,m-1) \times  
\prod_{X\in B(n)_m} {\cal U}{\cal D}(X)$$
have the same cardinality. 
\et

For $m=1$, Theorem \ref{jg}
gives $n|T(n,1)|=n^n$,  a result for which there is a celebrated bijective
proof {\bf\cite{j}}.

\noi{\bf Problem} {\em Find bijective proofs of Theorems \ref{gg} and
\ref{jg}}.

Recently, a related open problem, that of finding a combinatorial proof of
the product formula for the complexity of the hypercube was solved in {\bf\cite{b}}.


\begin{thebibliography}{AAA}


\bibitem{bvp} C. Bachoc, F. Vallentin and A. Passuello,
\newblock{\it Bounds for projective codes from semidefinite programming},
\newblock  arXiv: 1205.6406 (2012).   

\bibitem{b} O. Bernardi,
\newblock { \it On the spanning trees of the hypercube and other products of
graphs,}
\newblock Electronic J. Comb. {\bf 19(4)}, Paper 51 (16 Pages) (2012).



\bibitem{bh} A. E. Brouwer, and W. H. Haemers,
\newblock {\it Spectra of graphs},
\newblock  Springer, 2012.


\bibitem{cst} T. Ceccherini-Silberstein, F. Scarabotti, and F. Tolli,
\newblock {\it Harmonic analysis on finite groups},
\newblock  Cambridge University Press, 2008.



\bibitem{d1} P. Delsarte,
\newblock{\it Association schemes and $t$-designs in regular semilattices},
\newblock  J. Combinatorial Theory, Series A, 20 : 230-243 (1976).   


\bibitem{d2} P. Delsarte,
\newblock{\it Hahn polynomials, 
discrete harmonics, and $t$-designs},
\newblock SIAM J. Applied Math., 34 :  157-166 (1978).


\bibitem{du} C. F. Dunkl, 
\newblock{\it An addition theorem for 
some $q$-Hahn polynomials},
\newblock Monatsh. Math., 85 : 5-37 (1978).    


%\bibitem{gst} D. Gijswijt, A. Schrijver and H. Tanaka, 
%\newblock {\it New upper
%bounds for nonbinary codes based on the Terwilliger
%algebra and semidefinite programming},
%\newblock J. Comb. Theory, Ser. A, 113 :  1719-1731 (2006).



%\bibitem{g} J. T. Go, 
%\newblock {\it The Terwilliger algebra of the hypercube},
%\newblock Eur. J. Comb.,  23 : 399-429 (2002).



\bibitem{gr} J. Goldman and G. -C. Rota,
\newblock{ \it The number of subspaces of a vector space},
\newblock  in {\em Recent progress in Combinatorics} (Proc. Third Waterloo Conf.
on Combinatorics 1968), Academic Press : 75-83  (1969).    

 
%\bibitem{gll} R. L Graham, S. Y. R. Li and W. W. Li,
%\newblock{ \it On the structure of $t$-designs},
%\newblock SIAM J. Alg. Discr. Methods, 1 : 8-14 (1980).



%\bibitem{gj} J. E. Graver and W. B. Jurkat,
%\newblock{ \it The module structure of integral designs},
%\newblock J. Comb. Theory, Ser. A, 15 : 75-90 (1973).

\bibitem{hh} S. Hitzemann and W. Hochst\"{a}ttler,
\newblock{ \it On the combinatorics of Galois numbers},
\newblock  Discrete Math., 310 : 3551-3557 (2010).   

\bibitem{jl} G. James, and M.Liebeck,
\newblock {\it Representations and Characters of Groups},
\newblock  Cambridge University Press, 2001.

\bibitem{j} A. Joyal,
\newblock{ \it Une th\'{e}orie combinatoire des s\'{e}ries formelles},
\newblock  Advances in Math., 42 : 1-82 (1981).   



\bibitem{kc} V. Kac and P. Cheung,
\newblock {\it Quantum Calculus},
\newblock  Springer-Verlag, 2002.



%\bibitem{k} D. E. Knuth,
%\newblock{ \it Combinatorial Matrices},
%\newblock in {\em Selected papers on discrete mathematics}, 
%{\em CSLI lecture Notes, 106}, CSLI Publications, Stanford, CA: 177-186 (2003).

%\bibitem{mp1} J. M. Marco and J. Parcet,
%\newblock{ \it On the natural representation
%of $S(\Omega)$ into $L^2({\cal
%P}(\Omega))$: discrete harmonics and Fourier transform}, 
%\newblock J. Comb. Theory, Ser. A, 100 : 153-175 (2002).



\bibitem{mp2} J. M. Marco and J. Parcet,
\newblock{\it Laplacian operators and 
Radon transforms on Grassmann graphs},
\newblock Monatsh. Math., 150 : 97-132 (2007).

\bibitem{nsw} A. Nijenhuis, A. E. Solow and H. S. Wilf,
\newblock{\it Bijective methods in the theory of finite vector spaces},
\newblock J. Comb. Theory, Ser. A, 37 : 80-84 (1984).


\bibitem{p} R. A. Proctor,
\newblock{\it Representations of ${\mathfrak{sl}}(2,\C)$ 
on posets and the Sperner property},
\newblock SIAM J. Alg. Discr. Methods, 3 : 275-280 (1982).




\bibitem{s} A. Schrijver,
\newblock{\it New code upper bounds from the 
Terwilliger algebra and semidefinite programming},
\newblock IEEE Tran. Information Theory, 51 : 2859-2866 (2005).

%\bibitem{n} N. Singhi,
%\newblock{\it Tags on subsets},
%\newblock  Discrete Math., 306: 1610-1623 (2006).   

\bibitem{sr1} M. K. Srinivasan, 
\newblock{\it Symmetric chains, Gelfand-Tsetlin chains, 
and the Terwilliger algebra of the binary Hamming scheme},
\newblock J. Algebraic Comb., 34 : 301-322 (2011).

\bibitem{sr2} M. K. Srinivasan, 
\newblock{\it A positive combinatorial 
formula for the complexity of the $q$-analog of the $n$-cube},
\newblock Electronic J. Comb., 19(2) : Paper 34 (14 Pages) (2012).


\bibitem{sr3} M. K. Srinivasan,
\newblock{\it Notes on explicit block diagonalization},
\newblock  in {\em Combinatorial
Matrix Theory and Generalized Inverses of Matrices}, Springer: 13-31  (2013).


\bibitem{st1} R. P. Stanley,
\newblock{ \it Variations on differential posets},
\newblock in {\em Invariant Theory and Tableaux}, volume 19 of
{\em IMA Vol. Math. Appl.}, Springer: 145-165 (1990).



\bibitem{st2} R. P. Stanley,
\newblock {\it Enumerative Combinatorics - Volume 1, Second Edition},
\newblock  Cambridge University Press, 2012.




\bibitem{t} P. Terwilliger,
\newblock{ \it The incidence algebra of a uniform poset},
\newblock in {\em Coding theory and design theory, Part I}, volume 20 of
{\em IMA Vol. Math. Appl.}, Springer: 193-212 (1990).



%\bibitem{v} F. Vallentin, 
%\newblock {\it Symmetry in semidefinite programs},
%\newblock Linear Algebra Appl., 430: 360-369 (2009).


\end{thebibliography}
\end{document}